\documentclass[10pt,leqno]{amsart}
\usepackage{graphicx}
\usepackage{indentfirst, amssymb,amsmath,amsthm}
\usepackage{csquotes}
\usepackage{graphicx}
\numberwithin{equation}{section}
\begin{document}
\title[Sums of Powers of Natural numbers]{Proof Without Words: Sums of Powers of Natural numbers}
\date{}
\author[B. Chakraborty ]{Bikash Chakraborty }
\date{}
\address{Department of Mathematics, Ramakrishna Mission Vivekananda
Centenary College, Rahara, West Bengal 700 118, India. }
\email{bikashchakraborty.math@yahoo.com, bikashchakrabortyy@gmail.com}
\maketitle
\let\thefootnote\relax
\begin{abstract}
The aim of this short note is to give  world less proofs of the sum formula of $1^k+2^k+\ldots+n^k$, $k\in\{1,2,3\}$. Though the identity $1^2+2^2+...+n^2=\frac{n(n+1)(2n+1)}{6}$ is already published (B. Chakraborty, Proof without words: the sum of squares, Math. Intelligencer, 40 (2018), no. 2, pp. 20), yet for the sake of continuity of the sequel, we provide the visualization here again.
\end{abstract}
\section{$1+2+\ldots+n=\frac{n(n+1)}{2}$}
\begin{figure}[h!]
    \centering
    \includegraphics[scale=.8]{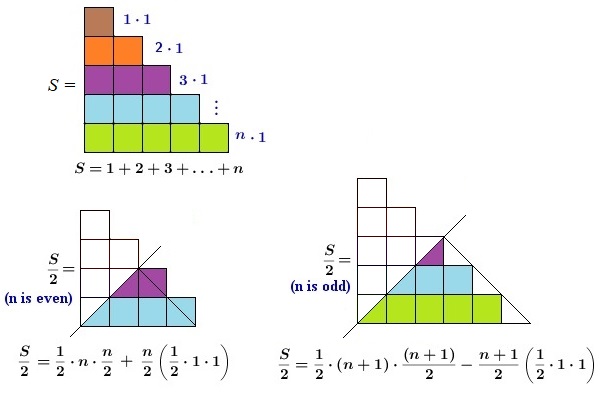}
    \end{figure}
\newpage
\section{$1^2+2^2+...+n^2=\frac{n(n+1)(2n+1)}{6}$}
\begin{figure}[h!]
    \centering
    \includegraphics[scale=.1]{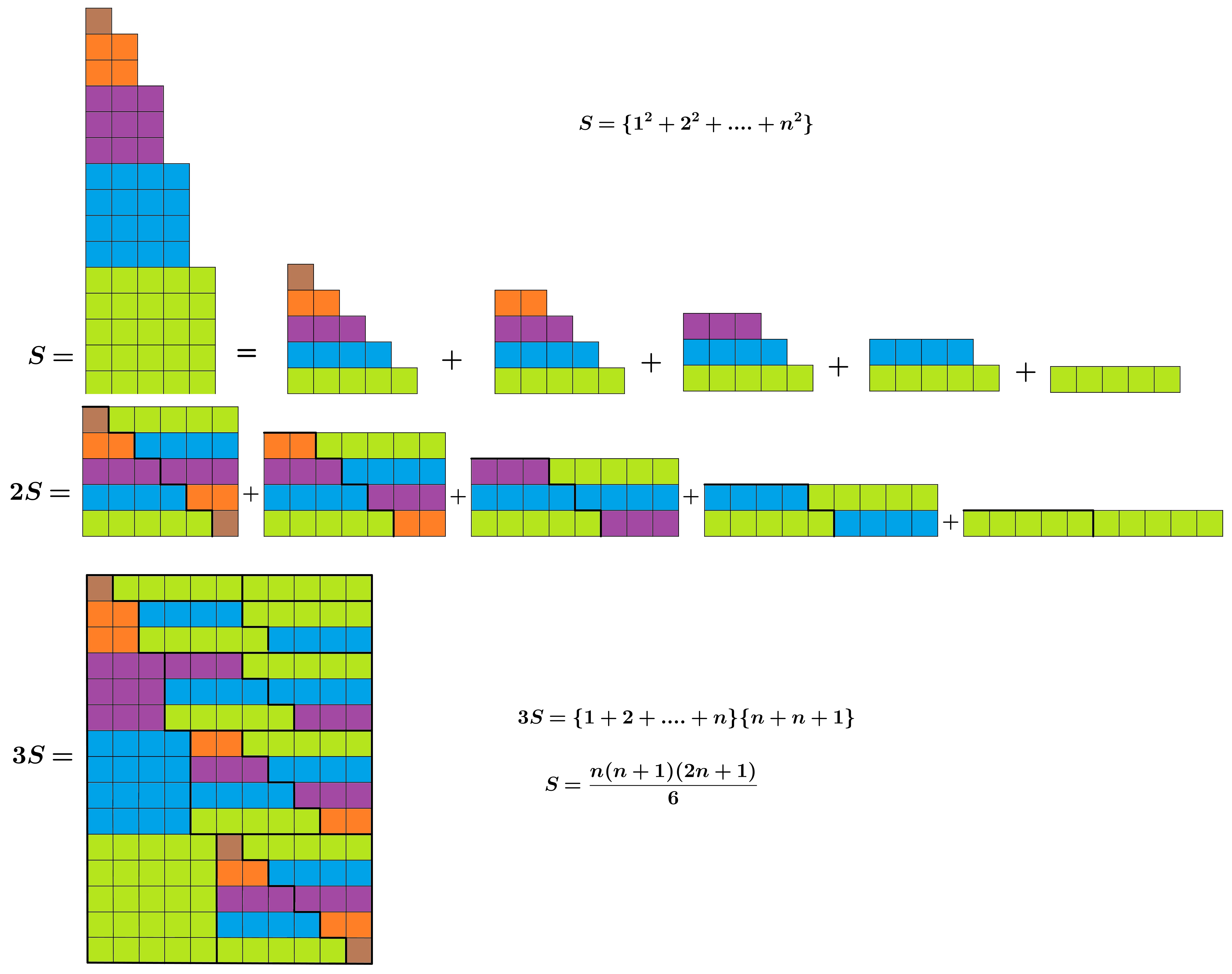}
    \end{figure}
    \footnote{The bibliographic data of the above visual proof is \enquote{B. Chakraborty, Proof without words: the sum of squares, Math. Intelligencer, 40 (2018), no. 2, pp. 20.}}
\newpage
\section{$1^3+2^3+...+n^3=\left(\frac{n(n+1)}{2}\right)^{2}$}
\begin{figure}[h!]
    \centering
    \includegraphics[scale=.4]{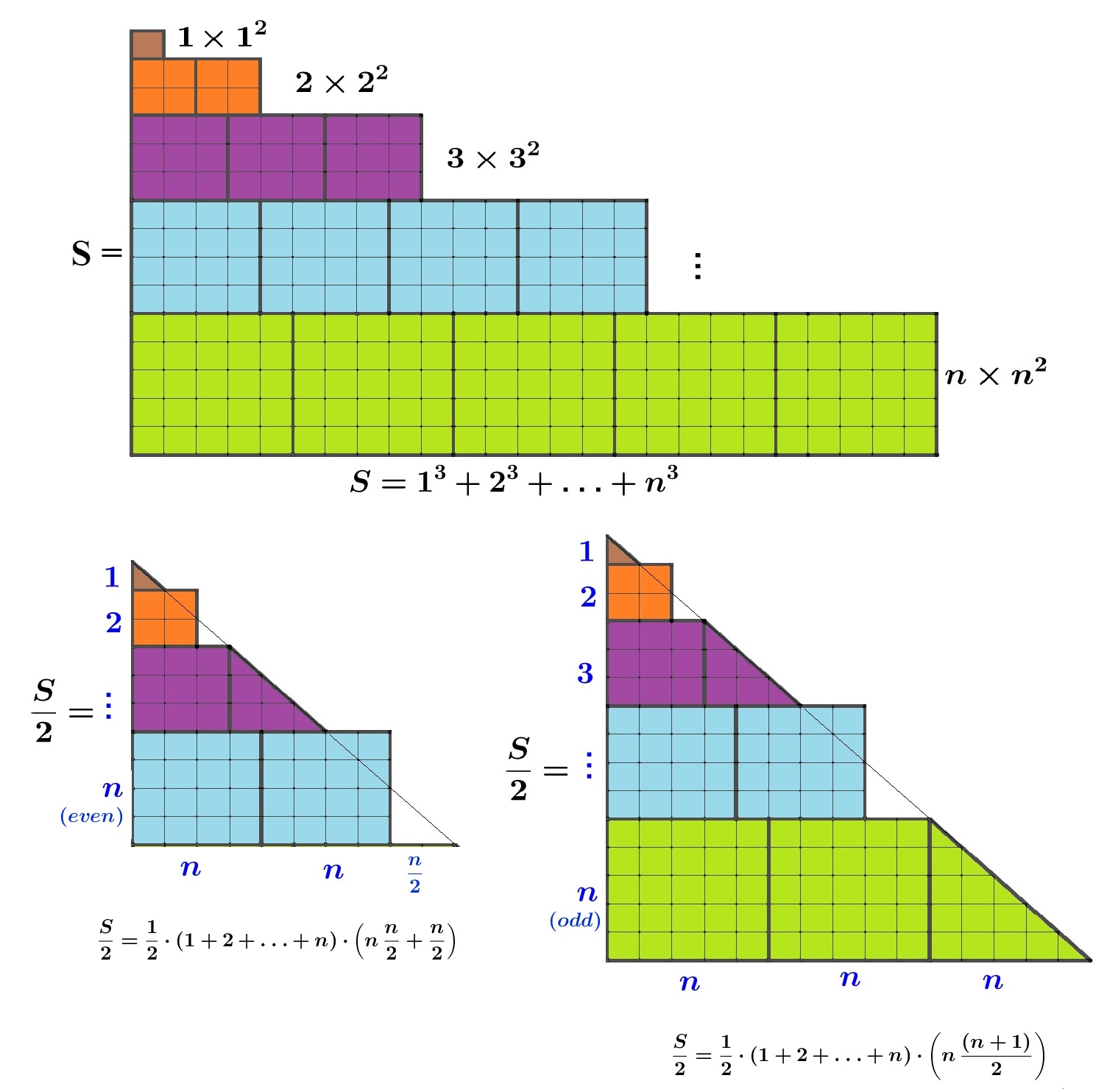}
    \end{figure}

\end{document}